\documentclass[11pt]{article}
\usepackage{amssymb,latexsym,amsmath}
\usepackage[english]{babel}
\usepackage{xcolor}
\usepackage{bm}
\usepackage{amsmath,amscd}
\usepackage{graphicx}
\usepackage{newtxtext}
\usepackage{newtxmath}
\usepackage{geometry}
\usepackage{hyperref}
\hypersetup{
    colorlinks = true,
    urlcolor   = blue,
    citecolor  = black,
}

\newcommand{\RomanNumeralCaps}[1]
\linenumbers

\title{Cosserat media in dynamics}

\author{\textbf{G. de Saxc\'e} \\
Univ. Lille, CNRS, Centrale Lille, UMR 9013 – LaMcube \\
Laboratoire de m\'ecanique multiphysique multi\'echelle, \\
F-59000, Lille, France, Email: gery.de-saxce@univ-lille.fr}

\begin{document}
\maketitle

\begin{abstract}
Our aim is to develop a general approach for the dynamics of material bodies of dimension $d$ represented by a mater manifold $\mathcal{N}$ of dimension $(d + 1)$ embedded into the space-time $\mathcal{M}$. It can be declined for $d = 0$ (pointwise object),  $d = 1$ (arch if it is a solid, flow in a pipe or jet if it is a fluid),  $d = 2$ (plate or shell if it is a solid, sheet of fluid), $d = 3$ (bulky bodies). We call torsor a skew-symmetric bilinear map on the vector space of affine real functions on the affine tangent space to the space-time. 
We use the affine connections as originally developed by Élie Cartan, that is the connections associated to the affine group. We introduce a general principle of covariant divergence free torsor from which we deduce 10 balance equations. We show the relevance of this general principle by applying it for $d$ from $1$ to $4$ in the context of the Galilean relativity.
\end{abstract}           


\vspace{0.2cm}








\section{Introduction}


Since a long time, scientists have been in search of tools to unify and structure the Mechanics, The most popular one today is probably the method of virtual powers or virtual works, initiated by Joseph-Louis Lagrange \cite{Lagrange 1788}, and the variational techniques. Without denying the power of these tools, our approach is different, using mainly the geometric methods. Our goal in this work is to revisit --for dynamical purposes-- the theory of generalized continua of the brothers Eugène and François Cosserat \cite{Cosserat, Cosserat 1909}.

Among the sources of inspiration of this work, one can quote Jean-Marie Souriau who highlighted the affine character of many features of the Mechanics and proposed a general approach called "Affine Mechanics" and based on tensor-distributions \cite{Souriau 1997a}. Starting from this approach, the attention is payed to the concept of affine tensor of which the importance in Mechanics was originally pointed out and developed by Wlodzimierz Tulczyjew 
and his school \cite{Tulczyjew 1988, Grabowska 2004}. In the present work, the balance equations are not obtained from the virtual power method but from the concept of parallel transport, using the affine connections as originally developed by Élie Cartan, that is the connections associated to the affine group \cite{Cartan 1923, Cartan 1924a}.

The paper is structured as follows.
\begin{itemize}
    \item In Section \ref{Section Affine tensors and torsors}, we define the torsor as a real  or vector-valued skew-symmetric bilinear map on the vector space of affine real functions on the affine tangent space to the space-time, a manifold of dimension 4. We give the basic notions of affine tensor calculus useful in the sequel, in particular the transformation law of the torsor components.
    \item The aim of Section \ref{Section Pointwise object} is to discuss the physical content of the torsors in the simplest case, the one of a pointwise object. We introduce the concept of Galilean connections and proper frames.
    \item In Section \ref{Section Cosserat media of arbitrary dimensions}, the previous approach is generalized  to Cosserat media of arbitrary dimension. For this, we need to consider a matter manifold of dimension  $(d + 1)$ . Before to go further, we verify that we recover classical results for the dynamics of a Cauchy medium. 
    \item To illustrate these general tools, to be more concrete, we particularize in Section \ref{Section 1D Cosserat media} to the Cosserat media of dimension $d = 1$ (beams and arches). We study the structure of the torsor of force-mass and we idealize 3D slender bodies thanks to a 3D-to-1D reduction.
    \item In Section \ref{Section Affine tensor analysis and principle of divergence free torsor}, we browse through the affine tensor analysis and we introduce a principle of covariant divergence free torsor.
    \item This principle is abstract but covers a broad spectrum of applications to the Dynamics that are presented in Section \ref{Section Dynamics of a bulky body: the case of a Cauchy medium} ($d = 3$, Cauchy medium), Section \ref{Section Dynamics of a pointwise object} ($d = 0$,  pointwise object), Section \ref{Section 1D Cosserat media}  ($d = 1$, beams and arches), Section \ref{Section Dynamics of 2D Cosserat media} ($d = 2$, plate and shells) and Section \ref{Section Dynamics of a 3D Cosserat medium} ($d = 3$, Cosserat medium)
    \item Section \ref{Section Conclusions and perspectives} is devoted to conclusions and perspectives.
\end{itemize}

\section{Affine tensors and torsors}
\label{Section Affine tensors and torsors}

By torsor, we mean a mathematical object able to modelize the behaviour of material bodies. In a nutshell, we present the general setting. For the Dynamics, the space-time $\mathcal{M}$  is the appropriate framework in Relativity but also in classical Mechanics. It is a manifold. $AT_{\bm{X}}\mathcal{M}$ is the affine space associated to the tangent vector space to $\mathcal{M}$ at $X$. $A^* T_{\bm{X}}\mathcal{M}$ is the vector space of the real-valued affine functions $\bm{\Psi}$ on $AT_{\bm{X}}\mathcal{M}$. They are called affine forms. A \textbf{torsor} is a bilinear and skew-symmetric function on this space
$$\bm{\tau} (\bm{\Psi}, \hat{\bm{\Psi}} ) = - \bm{\tau} (\hat{\bm{\Psi}} ,\bm{\Psi})
$$
It is real or vector-valued. This abstract definition may appear rather arbitrary but it will be justified in the sequel by various applications to the continuum mechanics. 

\subsection{Affine tensors}

The aim of this section is to give the basic notions of affine tensor calculus useful in the sequel. To know more about affine tensors and bundles, the reader referred to the publications of Tulczyjew, Urbañski and Grabowski \cite{Tulczyjew 1988, Grabowska 2004} and to author's one \cite{de Saxce 2003, de Saxce 2011, AffineMechBook}.

By the choice of an affine frame, that is an origin $\bm{a}_0$ and a basis $(\vec{\bm{e}}_\alpha)$ of $AT_{\bm{X}}\mathcal{M}$, the affine form  $\bm{\Psi}$ is characterized in a unique way by its height $\chi$ at the origin and the associated linear form $\bm{\Phi}$. The simplest affine tensors are:
\begin{itemize}
    \item Firstly, the \textbf{points} $\bm{a}$  of the affine space. Similarly to vectors, they are 1-contravariant. By the choice of an affine frame,  decomposing the bound vector from the origin $\bm{a}_0$ to $\bm{a}$ in the basis
    $$\bm{a} = \bm{a}_0 + V^\alpha \vec{\bm{e}}_\alpha
    $$
    we assign to $\bm{a}$ components $V^\alpha$.
    \item Secondly, the \textbf{affine forms} $\bm{\Psi}$. Similarly to linear forms, they are 1-covariant. Decomposing it in the basis of $A^* T_{\bm{X}}\mathcal{M}$, that is the constant function of value equal to 1, denoted $\bm{1}$, and the co-basis $(\bm{e}^\alpha)$  
    $$\bm{\Psi} = \chi \bm{1} + \Phi_\alpha \bm{e}^\alpha
    $$ 
    we assign to $\bm{\Psi}$ as components the height $\chi$ at the origin and the components $\Phi_\alpha$ of the unique associated linear form $\bm{\Phi}$.
\end{itemize}

Now, we present the affine tensors that will turn out to be the most relevant for the Mechanics.  
\begin{itemize}
\item The \textbf{torsors} that are 2-contravariant.
To assign components $(T^{\alpha}, J^{\alpha\beta})$ to $\bm{\tau}$, we decompose its two arguments $\bm{\Psi}$ and $\hat{\bm{\Psi}}$ in the basis of $A^* T_{\bm{X}}\mathcal{M}$
$$ \bm{\tau} (\bm{\Psi}, \hat{\bm{\Psi}} ) = 
\bm{\tau} (\chi \bm{1} + \Phi_\alpha \bm{e}^\alpha, \hat{\chi} \bm{1} + \hat{\Phi}_\beta \bm{e}^\beta ) 
$$
As usual, we use the bilinearity
\begin{eqnarray}
     \bm{\tau} (\bm{\Psi}, \hat{\bm{\Psi}} ) = 
 \chi  \hat{\chi}\, \bm{\tau}  (\bm{1}, \bm{1}) 
      +  \chi \, \hat{\Phi}_\beta \bm{\tau} (\bm{1} , \bm{e}^\beta)
      + \Phi_\alpha \hat{\chi} \, \bm{\tau} (\bm{e}^\alpha, \bm{1})
      + \Phi_\alpha \hat{\Phi}_\beta \bm{\tau} ( \bm{e}^\alpha, \bm{e}^\beta)
\label{tau(Psi, hat(Psi) =}
\end{eqnarray}
and the skew-symmetry then 
$$T^{00} = \bm{\tau}  (\bm{1}, \bm{1}) = 0, \;
              T^{\beta} = \bm{\tau}  (\bm{1},\bm{e}^\beta) = - \bm{\tau}  (\bm{e}^\beta, \bm{1}), \;
              J^{\alpha\beta} = \bm{\tau}  (\bm{e}^\alpha, \bm{e}^\beta) = - J^{\beta\alpha}$$
Moreover pay attention to the following subtleties. $\hat{\Phi}_\beta$
               is the value of $\bm{\hat{\Phi}}$  for the basis vector $\vec{\bm{e}}_\beta$ but, as $\bm{\hat{\Phi}}$ is the linear part of $\bm{\hat{\Psi}}$, by convention it is the value of $\bm{\hat{\Psi}}$ for $\vec{\bm{e}}_\beta$ and, identifying $\vec{\bm{e}}_\beta$ to an element of the bidual, it is the value of $\vec{\bm{e}}_\beta$ for $\bm{\hat{\Psi}}$ 
$$\hat{\bm{\Phi}}_\beta  = \hat{\bm{\Phi}}  (\vec{\bm{e}}_\beta) 
        = \hat{\bm{\Psi}} (\vec{\bm{e}}_\beta) =  \vec{\bm{e}}_\beta (\hat{\bm{\Psi}})$$
Similarly, identifying the origin $\bm{a}_0$ to an element of the bidual, $\chi$ is the value of $\bm{a}_0$ for $\bm{\Psi}$
$$\chi =\bm{\Psi}  (\bm{a}_0) = \bm{a}_0 (\bm{\Psi}) $$
then 
$$\chi \, \hat{\Phi}_\beta  = 
\bm{a}_0 (\bm{\Psi}  ) \, \vec{\bm{e}}_\beta (\hat{\bm{\Psi}}) 
= (\bm{a} _0 \otimes \vec{\bm{e}}_\beta)  (\bm{\Psi}, \hat{\bm{\Psi}} ) 
$$
Equation (\ref{tau(Psi, hat(Psi) =}) becomes
$$ \bm{\tau} (\bm{\Psi}, \hat{\bm{\Psi}} ) = 
      T^{\beta} ( \bm{a}_0 \otimes \vec{\bm{e}}_\beta - \vec{\bm{e}}_\beta  \otimes\bm{a}_0)  (\bm{\Psi}, \hat{\bm{\Psi}} )
      + J^{\alpha\beta} (\vec{\bm{e}}_\alpha \otimes \vec{\bm{e}}_\beta)  (\bm{\Psi}, \hat{\bm{\Psi}} )
$$
As $\bm{\Psi}$ and $\hat{\bm{\Psi}}$ are arbitrary,
we obtain the decomposition 
$$\bm{\tau} = 
      T^{\beta}\, ( \bm{a}_0 \otimes \vec{\bm{e}}_\beta - \vec{\bm{e}}_\beta  \otimes\bm{a}_0)
     +  J^{\alpha\beta}\,\vec{\bm{e}}_\alpha  \otimes \vec{\bm{e}}_\beta
     $$
\item The \textbf{co-torsors} are 2-covariant and can be put in duality with the torsors (for more details, see \cite{AffineMechBook}, Section 5.1.3., page 76-80).
\item The \textbf{momentum tensors} are mixed tensors, 1-covariant and 1-contravariant (see \cite{AffineMechBook}, Section 16.3., page 325-327). 
\end{itemize} 
These two latter kinds of affine tensors are also useful in Mechanics and we shall give some words at the end of the paper.

\subsection{Transformation laws of affine tensors}
\label{Subsection Transformation laws of affine tensors}

Now, we discuss the transformation laws of affine tensors, considering a change of affine frames  $(\bm{a}_0, (\vec{\bm{e}}_\alpha)) \longrightarrow (\bm{a}'_0, (\vec{\bm{e}}'_\beta)) $ given by the components $C'^\beta$ of the bound vector pointing from $\bm{a}'_0$ to $\bm{a}_0$  in the new basis, that can be stored in a column $C'$, and the transformation matrix $P$ of the basis change
\begin{eqnarray}
    \overrightarrow{\bm{a}'_0 \bm{a}_0} = C'^\beta \vec{\bm{e}}'_\beta, \qquad 
  \vec{\bm{e}}'_\beta = P^\alpha_\beta \vec{\bm{e}}_\alpha
\label{affine frame change}
\end{eqnarray}
Taking into account that the space-time is of dimension $4$, we have in matrix form
\begin{eqnarray}
     C' =\left( \begin{array} {c}
                       C'^1    \\
                       \vdots\\
                       C'^4   \\
                    \end {array} \right), \qquad 
   P =\left( {{\begin{array}{*{20}c}
                 P^1_1 \hfill & \ldots \hfill & P^1_4 \hfill \\
                 \vdots \hfill & \ddots \hfill & \vdots \hfill \\
                 P^4_1 \hfill & \ldots \hfill & P^4_4 \hfill \\
              \end{array} }} \right), \qquad
              C = -P \, C'
\label{C' = & P = & C =}
\end{eqnarray}
Then —in the matrix form— the laws of transformation of affine tensors are given by

\begin{itemize}
    \item for a point $\bm{a}$ of components $V^\alpha$ by
    $$V' = C' + P^{-1} V
    $$
    storing its components into the column $V$,
    \item for an affine form $\bm{\Psi}$ of components $(\chi, \Phi_\alpha)$ by
    $$\chi' = \chi + \Phi\,C,\qquad \Phi' = \Phi\,P 
    $$
    gathering the components $\Phi_\alpha$ into the row $\Phi$,
    \item for a torsor $\bm{\tau}$ of components $(T^\alpha,  J^{\alpha\beta})$ by
    $$ T' =P^{-1} T, \qquad J' = P^{-1}  (J + C T^T - T C^T) \, P^{-T}  
    $$
    storing the components $T^\alpha$ into the column $T$ and $J^{\alpha\beta}$ into the matrix $J$. It is worth to remark that the components $T^\alpha$ represent a vector. It is the linear part of the torsor.
\end{itemize}
\vspace{0.3cm}
Now, combining the change of affine frames, we observe that the set of couples $(C, P)$ is the affine group $\mathbb{GA} (4) = \mathbb{R}^4 \rtimes \mathbb{GL} (4)$. For a more rational organization of calculations, it is convenient to use a classical trick, the  linear representation on $\mathbb{R}^5$ of the affine group of $\mathbb{R}^4$
$$     \tilde{P} =\left( {{\begin{array}{*{20}c}
                 1 \hfill & 0 \hfill \\
                 C \hfill & P \hfill \\
              \end{array} }} \right), \qquad
$$
With this extra fifth dimension without physical meaning, we store the components of affine tensors into bigger column, row and matrix 
$$ \tilde{V} =\left( \begin{array} {c}
                       1    \\
                       V    \\
                    \end {array} \right), \qquad 
 \tilde{\Psi} = (\chi, \; \Phi), \qquad 
 \tilde{\tau}     =\left( {{\begin{array}{cc}
                          0  &  T^T \\
                        - T  &  J   \\
                      \end{array} }} \right)
$$

Then the transformation laws of the corresponding tensors take these simpler and compact forms
$$ \tilde{V} '= \tilde{P}^{-1} \tilde{V}, \qquad
      \tilde{\Psi}' = \tilde{\Psi}\,\tilde{P}, \qquad
      \tilde{\tau}' = \tilde{P}^{-1} \tilde{\tau} \, \tilde{P}^{-T}
$$

\section{Pointwise object}
\label{Section Pointwise object}

Our aim now is to discuss the physical content of the torsors in the simplest case, the one of a pointwise object that may be a material particle or an object which has an internal structure but can be thought as pointwise if it is seen from a long way off (for instance a satellite seen from the Earth).

\subsection{Galilean mechanics}

The element of the space-time, the events are represented in a suitable local chart by these coordinates, the time $t$ and the position $x$, 
stored in the 4-column 
$$ X =\left( \begin{array} {c}
                       t \\
                       x \\
                    \end {array} \right) \in \mathbb{R}^4
                    $$
The Galilean transformations are the transformations preserving inertial motions, durations, distances and volume orientations, 
then affine of the form  $X  = P\,X' + C$  with
$$
P =\left( \begin{array} {cc}
                        1  & 0 \\
                        u & R
             \end {array} \right),\qquad
C =\left( \begin{array} {c}
                        \tau_0 \\
                        k
             \end {array} \right)         
$$
where $P$ is a 4-by-4 matrix, $u$ the Galilean boost or velocity of transport, $R$ a rotation and $C$ a space-time translation. Their set is the Galilei group, a Lie group of dimension $10$. By the way, it may be wondered why we specify the variance of tensors. It is because the Galilean geometry is not Riemannian. We are not allowed to lower or raise the indices. 

\subsection{Proper frame}

\begin{figure}[h!]
\centering
\includegraphics[scale=.80]{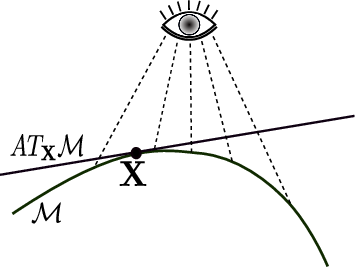}
\caption{The manifold perceived as an affine space}
\label{fig Manifold as affine space}
\end{figure}

As said Élie Cartan \cite{Cartan 1923, Cartan 1924a} (Figure \ref{fig Manifold as affine space}): \\
“The affine space at $\bm{X}$ could be seen as the manifold itself
that would be perceived in an affine manner by an observer located at $\bm{X}$”

The zero tangent vector $\bm{0}$ at $\bm{X}$ --that is the equivalence class of curves reduced to $\bm{X}$-- can be identified to $\bm{X}$ itself. It is a particular origin. A frame of reference in which an object is \textbf{at rest} at the special origin $\bm{0}$  of the affine space $AT_{\bm{X}}\mathcal{M}$ is called a \textbf{proper frame} attached to this object. 

\subsection{The torsor of a pointwise object}

\label{Subsection The torsor of a pointwise object}

The torsor of a pointwise object is scalar-valued.

For an object in a proper frame, we claim that $T'$ and $J'$ have a reduced form with the only non vanishing components $m$, a scalar, and the axial vector $l_0$ associated to the skew-symmetric matrix $j(l_0)$ 
$$ \tilde{\tau}'     =\left( {{\begin{array}{cc}
                          0  &  T'^T \\
                        - T'  &  J'   \\
                      \end{array} }} \right)
                = \left( {{\begin{array}{ccc}
                          0     & m     &   0      \\
                        - m     &  0     &   0       \\
                          0     &  0     & - j (l_0) \\
                      \end{array} }} \right)
$$

Changing of affine frame, we apply a Galilean boost of velocity $v$ and a translation of the origin at $x$ 
\begin{eqnarray}
     \tilde{P}  = \left( {{\begin{array}{cc}
                          1  &  0 \\
                          C  &  P \\
                      \end{array} }} \right)
                  = \left( {{\begin{array}{ccc}
                          1     &  0     &   0      \\
                          0     &  1     &   0      \\
                          x     &  v     & 1_{\mathbb{R}^3} \\
                      \end{array} }} \right)
\label{tilde(P) = ((1 0)(C P)) =}
\end{eqnarray}
the transformation law of torsors gives in the new frame of reference this representative 
$$ \tilde{\tau}   = \tilde{P}  \, \tilde{\tau}' \tilde{P}^T
                       = \left( {{\begin{array}{ccc}
                          0     &  m      &   p^T      \\
                        - m      &  0     &   - q^T       \\
                        - p     &   q     & - j (l) \\
                      \end{array} }} \right)
$$
where $m$ is invariant, identified to the mass and there are new non-vanishing components:
\begin{itemize}
    \item the linear momentum $p = m\, v$,
    \item the quantity of position $q = m\, x$ (called passage by Souriau in \cite{SSD, SSDEng}), the forgotten momentum, unknown to standard textbooks,
    \item the angular momentum $l = l_0 + x \times m\,v$, decomposed into the proper angular momentum $l_0$ and the orbital angular momentum. 
\end{itemize}
\vspace{0.3cm}
It is worth to remark that the expression $m \, v$ of $p$ and the law of transport of the angular momentum are not definitions or postulates but by-products of the transformation law of affine tensors.

\section{Cosserat media of arbitrary dimensions}
\label{Section Cosserat media of arbitrary dimensions}

\begin{figure}[h!]
\centering
\includegraphics[scale=.60]{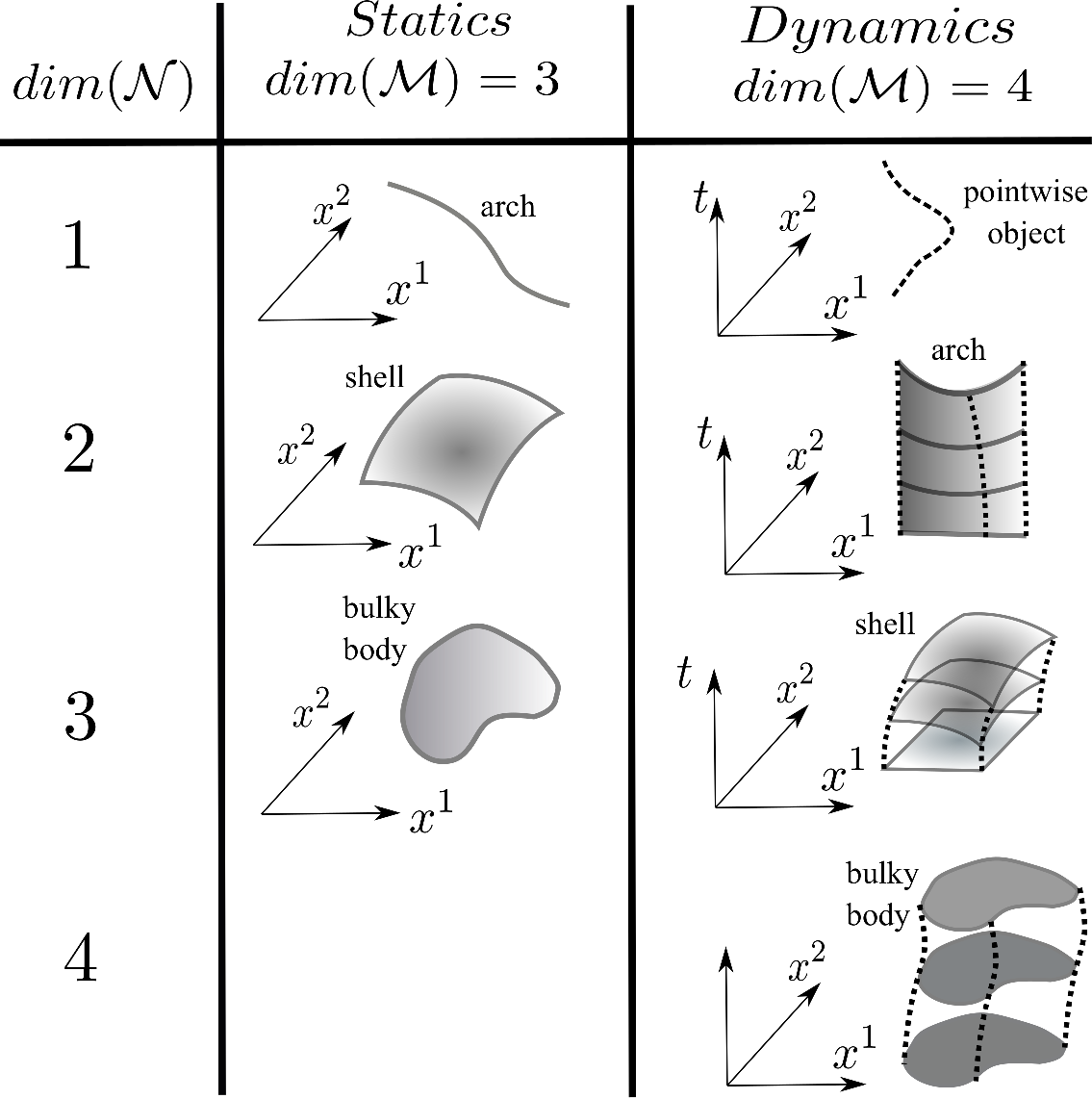}
\caption{Cosserat media of arbitrary dimensions}
\label{fig Material bodies}
\end{figure}

Now, we would like to generalize the previous approach to Cosserat media of arbitrary dimension. For this, we need to consider a \textbf{matter manifold} $\mathcal{N}$. Different applications are represented in Figure \ref{fig Material bodies}.

Before to consider the dynamics, let us have a glance at the limit case of Statics, cancelling the time, then $\mathcal{M}$ is reduced to the physical space of dimension 3. To modelize an arch —which is geometrically a curve—, the matter manifold $\mathcal{N}$ is of dimension 1, to modelize a shell, $\mathcal{N}$ is of dimension 2, and for the bulky body, of dimension 3.

Next let us come back to the Dynamics by restoring the time dimension. For a particle, the matter manifold $\mathcal{N}$ is the trajectory in the space-time, then of dimension 1. For an arch, $\mathcal{N}$ is of dimension 2. The trajectories of particles are dotted lines and the shapes of the arch at different times are solid lines. Below, $\mathcal{N}$ is of dimension 3 for a shell and dimension 4 for a bulky body.

\subsection{The torsor of a Cosserat medium of arbitrary dimension}

\begin{figure}[h!]
\centering
\includegraphics[scale=.60]{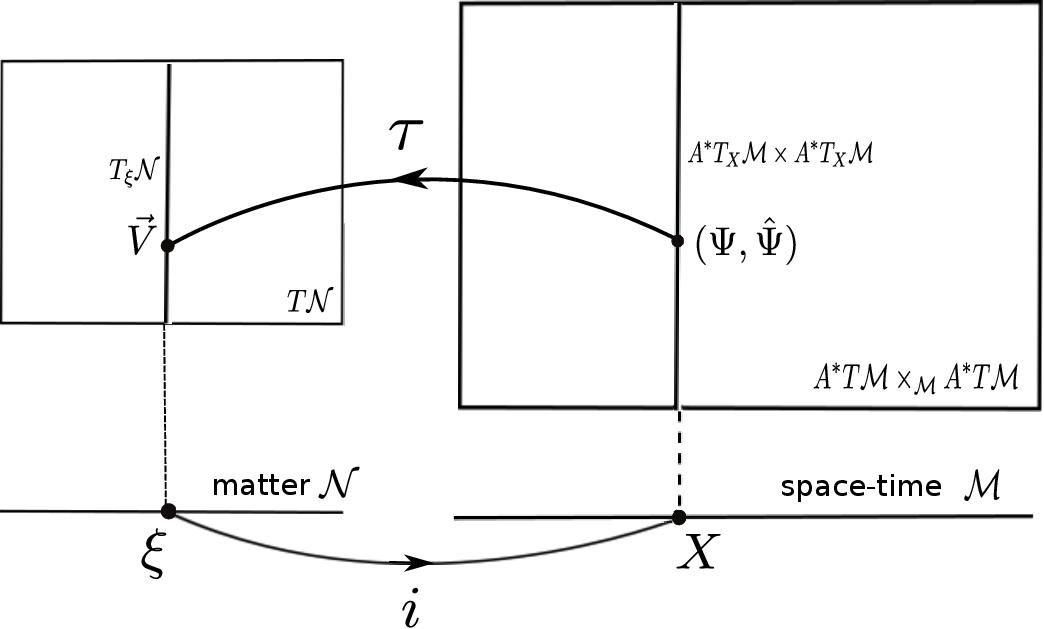}
\caption{Vector-valued torsor}
\label{fig Vector valued torsor}
\end{figure}

To generalize the torsor of a pointwise object to Cosserat media of arbitrary dimension, we must consider that it is vector-valued (Figure \ref{fig Vector valued torsor}). The matter manifold $\mathcal{N}$ is the set of material particles $\bm{\xi}$ and we describe the motion of the matter by an embedding $i$ from $\mathcal{N}$ into the space-time $\mathcal{M}$. We consider a torsor at $\bm{X} = i (\bm{\xi})$, with arguments $\bm{\Psi}, \bm{\hat{\Psi}}$,  valued in the tangent vector space to $\mathcal{N}$ at $\bm{\xi}$. Then, we build this pullback bundle of the product of $A^* T_{\bm{X}}\mathcal{M}$  by itself by moving the fiber over $\bm{X}$ to a fiber over $\bm{\xi}$, next we define the following bundle map over $\mathcal{N}$
$$ i^* (A^* T\mathcal{M} \times_{\mathcal{M}} A^* T\mathcal{M} ) \rightarrow T\mathcal{N} 
$$
By convention, the space-time indices (related to $\mathcal{M}$) are at the right and the material indices (related to $\mathcal{N}$) are at the left. In an affine frame $(\bm{a}_0, (\vec{\bm{e}}_\alpha))$  of $A^* T_{\bm{X}}\mathcal{M}$  and a basis $(_\gamma \vec{\bm{\eta}})$ of $T_{\bm{\xi}}\mathcal{N}$, the torsor $\bm{\tau}$ is decomposed as follows
$$ \bm{\tau} = \,^\gamma \bm{\tau}\,_\gamma\vec{\bm{\eta}},\qquad \,^\gamma \bm{\tau} =
      \,^{\gamma}T^{\beta}\, ( \bm{a}_0 \otimes \vec{\bm{e}}_\beta - \vec{\bm{e}}_\beta  \otimes\bm{a}_0)
     + \,^\gamma J^{\alpha\beta}\,\vec{\bm{e}}_\alpha  \otimes \vec{\bm{e}}_\beta
$$
where the value of $\bm{\tau}$ is given in the basis $(_\gamma \vec{\bm{\eta}})$ and the linear combination coefficient $^\gamma \bm{\tau}$ are decomposed into the affine frame as discussed before. Now, $T$ components have 2 indices and $J$ components have 3 indices.

\subsection{The torsor of a Cauchy medium}

At this stage, we want, before to go further, to verify that we recover classical results for the dynamics of a 3D body.
As the dimension of the matter manifold $\mathcal{N}$ and the space-time $\mathcal{N}$ are the same, we can choose for easiness local charts such that the embedding is represented by the identity
\begin{eqnarray}
    \xi^\alpha = X^\alpha
\label{xi^alpha = X^alpha}
\end{eqnarray}
then we can move the material index $\gamma$ at the right
\begin{eqnarray}
      ^{\gamma}T^{\beta} =T^{\beta\gamma}, \quad ^{\gamma}J^{\alpha\beta} =J^{\alpha\beta\gamma}
\label{convention index at the right}
\end{eqnarray}
In a proper frame attached to the elementary reference volume
\begin{itemize}
    \item if the $J^{\alpha\beta\gamma}$ components do not vanish, we modelize a full 3D Cosserat medium
    \item in contrast, if the $J^{\alpha\beta\gamma}$ components vanish, we modelize a 3D Cauchy medium. 
\end{itemize}
For a Cauchy medium in a proper frame, we claim that the symmetric 4-by-4 matrix $T'$ of which the components are the $T'^{\beta\gamma}$ is reduced to 
$$  T'  =  \left( {{\begin{array}{*{20}c}
                     \rho \hfill &   0       \hfill \\
                     0    \hfill & - \sigma \hfill \\
           \end{array} }} \right)
$$
where $\rho$ is the mass density and $\sigma$ is Cauchy's stress tensor. Applying a Galilean boost of velocity $v$, 
the transformation law of vector-valued torsors gives 
\begin{eqnarray}
     T  =  P \, T' P^T
      =  \left( {{\begin{array}{cc}
             \rho  &    p^T            \\
             p     &  - \sigma_{\star} \\
         \end{array} }} \right)
         =  \left( {{\begin{array}{cc}
             \rho    & \rho\,v^T              \\
             \rho\,v &  \rho\,v\,v^T - \sigma \\
         \end{array} }} \right)
\label{T = P T' P^T =}
\end{eqnarray}
where we recognize the linear momentum $p$ per volume unit
and the dynamical stresses $\sigma_{\star}$, well known in Fluid mechanics, then there is nothing new, it is just a control. For these reasons, the tensor of components $ T^{\beta\gamma}  $  is called  \textbf{stress-mass tensor}. 

\section{1D Cosserat media}
\label{Section 1D Cosserat media}

The tools presented here are general but, to illustrate, to be more concrete, we particularize now to the dimension 1. But what is meant by 1 dimensional Cosserat media? It may be a beam, an arch or a string if it is a solid, but it may be also a flow in a pipe, in a garden hose or a water jet, if it is a fluid. In the spirit of works such as \cite{Lighthill 1960, Boyer 2009, Boyer 2013} to represent the motion of a swimming fish, both the slender solid body and the surrounding flow can be considered. Our ambition is to develop a mathematical modelling covering all these applications. 

\subsection{Modelling of the matter motion}

The motion of a 1D material body can be described by the embedding $i$ of a matter manifold $\mathcal{N}$ of dimension 2 into the space-time $\mathcal{M}$. 

Let $s_0$  be the arclength with respect to a given reference point of the initial slender body (\textit{i.e.} at $t = 0$). In the Lagrangean coordinates $t$ and $s_0$ and the Galilean coordinates $X$, the position at time $t$ of the material particle pinpointed by $s_0$ is represented by
$$ x = \varphi (t, s_0)
$$
Its velocity is
\begin{equation}
   v = \frac{\partial \varphi}{\partial t}
\label{v = partial varphi / partial t} 
\end{equation}

\begin{figure}[h!]
\begin{center}
\includegraphics[scale=.60]{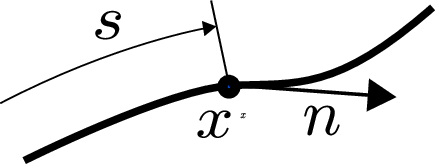} 
\caption{1D material body at time $t$}
\label{fig 1Dsxn}
\end{center}
\end{figure}

Let $s$ be the arclength with respect to a given reference point of the same slender body at time $t$ (Figure \ref{fig 1Dsxn}). In the Eulerian coordinates $t$ and $s$, the motion is represented by $x = \psi (t, s)$. We claim there exists a function 
\begin{equation}
   s_0 = f (t, s)
\label{s_0 = f (t, s)} 
\end{equation}
monotonically strictly increasing with respect to $s$ such that
\begin{equation}
   x = \psi (t, s) = \varphi (t, f (t, s))\ .
\label{x = psi (t, s) = varphi (t, f (t,s)} 
\end{equation}
Then the local representation of the embedding is
 $$ i: \xi =  \left( \begin{array} {c}
        t    \\
        s    \\
\end {array} \right) \mapsto X = 
 \left( \begin{array} {c}
        t    \\
        x    \\
\end {array} \right)  =
 \left( \begin{array} {c}
        t    \\
        \psi (t, s)   \\
\end {array} \right) 
 $$
Let us determine its tangent map $\bm{U}$. Differentiating (\ref{x = psi (t, s) = varphi (t, f (t,s)}) and taking into account (\ref{v = partial varphi / partial t}), it is represented by
$$\left( \begin{array} {c}
        dt    \\
        dx    \\
\end {array} \right) =
 \left( {{\begin{array}{*{20}c}
 1 \hfill                                                                          & 0 \hfill \\
 v + \dfrac{\partial f}{\partial t}\,\dfrac{\partial \varphi}{\partial s_0} \hfill &
     \dfrac{\partial f}{\partial s}\,\dfrac{\partial \varphi}{\partial s_0} \hfill            \\
\end{array} }} \right)\,
\left( \begin{array} {c}
        dt    \\
        ds    \\
\end {array} \right)
$$
As $s$ is the arclength, $ds = \parallel dx \parallel $ at constant $t$, hence
$$ n = \dfrac{\partial \psi}{\partial s}
     = \dfrac{\partial f}{\partial s}\,\dfrac{\partial \varphi}{\partial s_0}
$$
is the unit tangent vector to the curve. Differentiating (\ref{s_0 = f (t, s)}) at constant $s_0$ and owing to (\ref{v = partial varphi / partial t}), we obtain the tangential component of the velocity 
$$ v_t = v\cdot n = \frac{d s}{d t} \mid_{s_0 = C^{te}}
       = - \frac{\partial f}{\partial t}\ /\ \frac{\partial f}{\partial s}
$$
Thus, the tangent map $\bm{U}$ is represented by the 4-by-2 matrix
$$ dX = \left( \begin{array} {c}
        dt    \\
        dx    \\
\end {array} \right) =
  \left( {{\begin{array}{*{20}c}
                    1 \hfill          & 0 \hfill \\
                    v - v_t\,n \hfill & n \hfill \\
                \end{array} }} \right)\,
 \left( \begin{array} {c}
        dt    \\
        ds    \\
\end {array} \right) 
= U\,d\xi               
$$

\subsection{The tensor of force-mass}

The components $ ^{\gamma}T^{\beta}$ can be stored in a 2-by-4 matrix. In a proper frame attached to an elementary segment of  1D material body, it is reduced to
$$          T'  =  \left( {{\begin{array}{cc}
                      \rho_l &   0    \\
                      0      & - F^T \\
                  \end{array} }} \right)
$$
where $\rho_l$ is the mass density, $F$ is the statical force density acting through the cross-section (both per unit arclength). 

\begin{figure}[h!]
\begin{center}
\includegraphics[scale=.60]{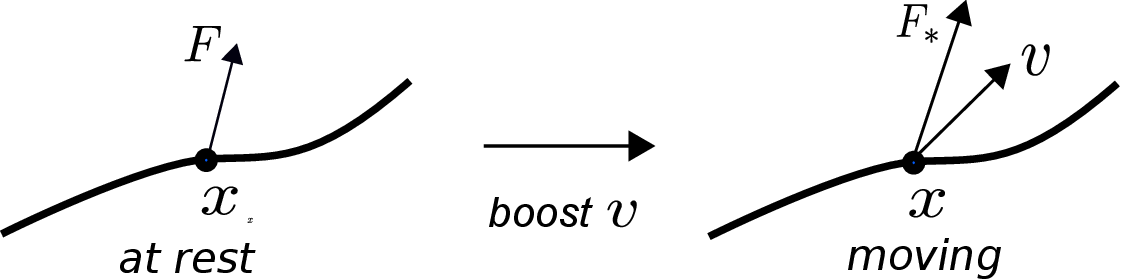} 
\caption{Boost $v$}
\label{fig 1Dboost-en}
\end{center}
\end{figure}

Next, applying a Galilean boost $v$ (Figure \ref{fig 1Dboost-en}), the transformation law of tensor gives 
$$ T  =  \left( {{\begin{array}{cc}
                  \rho_l  &    p^T         \\
                  p_t     &  - F^T_{\star} \\
            \end{array} }} \right)
         =  \left( {{\begin{array}{cc}
                  \rho_l     &    \rho_l v^T         \\
                  \rho_l v_t &  (\rho_l v_t v - F)^T \\
            \end{array} }} \right)
$$
with the new components of $T$:
\begin{itemize}
    \item $\rho_l$  which is invariant
    \item there are new non-vanishing components:
the linear momentum $p$ (per unit arclength)
and its tangential projection $p_t$ onto the curve
    \item finally the dynamical force $F_\star$ (something similar to the dynamical stress of the 3D body)
\end{itemize}
Naturally, $T$ represents a tensor called the force-mass tensor. 

\subsection{3D-to-1D reduction}
\label{Subsection 3D-to-1D reduction}

\begin{figure}[h!]
\begin{center}
\includegraphics[scale=.60]{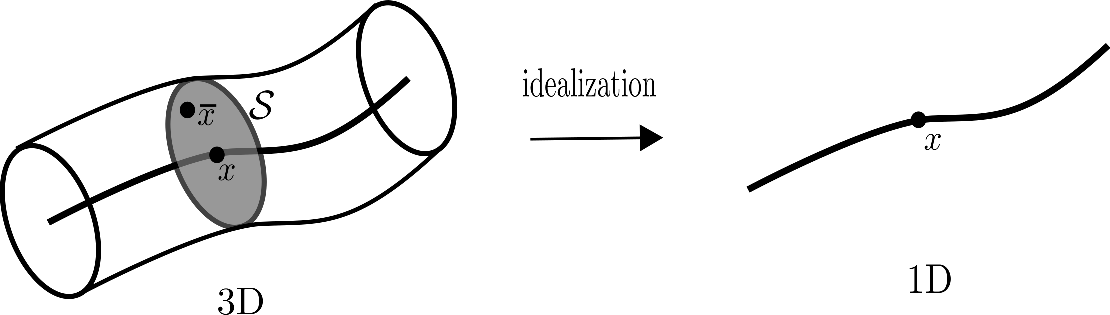} 
\caption{3D-to-1D reduction of a slender body}
\label{fig 3Dto1D-en}
\end{center}
\end{figure}

To modelize material bodies which are 3 dimensional by nature but slender with the size of the cross-section very small with respect to its length, engineers have idealized them since a long time --to shorten the calculations-- by 1 dimensional media, thanks to a 3D-to-1D reduction (Figure \ref{fig 3Dto1D-en}).

For this aim, we define a projection $\bm{\Pi}$ of $T_{\bm{X}} \mathcal{M}$ into $T_{\bm{\xi}} \mathcal{N}$ where $\bm{X} = i (\bm{\xi})$, given by a 2-by-4 matrix $\Pi$
$$  d\xi = \left( \begin{array} {c}
        dt    \\
        ds    \\
\end {array} \right) =
  \left( {{\begin{array}{*{20}c}
                    1 \hfill & 0   \hfill \\
                    0 \hfill & n^T \hfill \\
                \end{array} }} \right)\,
 \left( \begin{array} {c}
        dt    \\
        dx    \\
\end {array} \right)
= \Pi\,dX
$$
such that the product of the tangent map $U$ to the embedding and this projection $\Pi$ is the identity
$$ \Pi\,U = 1_{\mathbb{R}^2}
$$
In the sequel, quantities topped by a bar are related to the 3D body.
As discussed before, the 3D body, supposed to be a Cauchy medium, is described by a stress-mass tensor, owing to (\ref{T = P T' P^T =})
$$ \bar{T} =  \left( {{\begin{array}{cc}
             \rho          & \rho\,\bar{v}^T                   \\
             \rho\,\bar{v} & \rho\,\bar{v}\,\bar{v}^T - \sigma \\
         \end{array} }} \right)
$$

Combining $\bar{T}$ with the projection $\Pi$ and integrating on the cross-section $\mathcal{S}$, we obtain the components of the force-mass tensor $T$ of the 1D body
\begin{eqnarray}
    T = \int_{\mathcal{S}}\,\Pi \, \bar{T}\,d\mathcal{S}
\label{T = int_S Pi bar(T) dS}
\end{eqnarray}
where 
\begin{itemize}
    \item the component $\rho_l$  is, as expected, the mass per arclength
    $$\rho_l = \int_{\mathcal{S}}\, \rho\,d\mathcal{S}
    $$
    \item the linear momentum per arclength is given by
    $$ \rho_l v = \int_{\mathcal{S}}\,\rho\,\bar{v}\,d\mathcal{S}
    $$ 
    \item the force 
    \begin{eqnarray}
       F =  \int_{\mathcal{S}}\,(\sigma\,n - \rho\,(v_t - \bar{v}_t)(v - \bar{v}))\,d\mathcal{S}
       \label{F = int_S (sigma n - rho (v_t - bar(v)_t) (v - bar(v)) dS}
    \end{eqnarray}
    contains, besides the classical first term, a dynamical contribution resulting from the velocity fluctuations of $\bar{v}$ around the average value $v$ in the cross-section
\end{itemize}

\vspace{0.3cm}

Next we want to calculate the $J$ components by the reduction technique as for the $T$ ones.  
Our conventions are that the index of the time component is $0$, the indices of the spatial ones are $1, 2$  and $3$, the Greek indices run from 1 to 4 while the Latin  indices run from 1 to 3.
The torsor of the 3D body considered as a  Cauchy medium is such that the $\bar{J}$ components vanish at the points $\bar{x}$ of the cross-section, in the proper frame attached to the volume element at $\bar{x}$. Applying a spatial translation $\bar{x}$, the transformation law gives non vanishing $\bar{J}'$ components at the mass-center of the cross-section taken as new origin ($x = 0$) . Next, we skip the primes for easiness, that leads to
$$  \bar{J}^{i0\rho}  = \bar{x}^i \bar{T}^{0\rho},\qquad 
     \bar{J}^{ij\rho} = \bar{x}^i \bar{T}^{j\rho} 
     - \bar{x}^j \bar{T}^{i\rho}
$$
By a formula similar to (\ref{T = int_S Pi bar(T) dS}), we obtain the $J$ components 
\begin{eqnarray}
    ^\gamma J^{\alpha\beta} = \int_{\mathcal{S}}\,^\gamma \Pi_\rho \, \bar{J}^{\alpha\beta\rho}\,d\mathcal{S}
\label{J = int_S Pi bar(J) dS}
\end{eqnarray}
where we had to give up the matrix notation because the $J$ components have 3 indices. For the same reason, we reduce the number of indices by using engineer notations, putting
$$ q^i =\, ^0 J^{i0},\quad 
\quad l^i = \,^{0}J^{kl}, \quad
l^i_\star = \,^{1}J^{i0}, \quad
\quad M^i_\star = \,^{1}J^{kl} 
$$
where $(ikl)$ is a cyclic permutation of $(123)$,
and we obtain in matrix form for the 1D Cosserat medium
\begin{itemize}
    \item the position quantity
    $$ q = \int_{\mathcal{S}} \rho\,\bar{x}\,d\mathcal{S}
    $$
    \item the classical angular momentum
    $$  l = \int_{\mathcal{S}} \bar{x} \times \rho\, \bar{v}\,d\mathcal{S}
    $$
    \item a new momentum
    $$ l_\star = \int_{\mathcal{S}} \rho\,\bar{v}_t \bar{x}\,d\mathcal{S} 
    $$
    which is not obvious, unknown to standard textbooks
    \item and the dynamical moment
    $$ M_\star = \int_{\mathcal{S}} \bar{x} \times (\rho \bar{v}_t \bar{v} - (\sigma\,n))\,d\mathcal{S}
    $$
    analogous to the dynamical stress of the 3D body.
\end{itemize}

\begin{figure}[h!]
\begin{center}
\includegraphics[scale=.60]{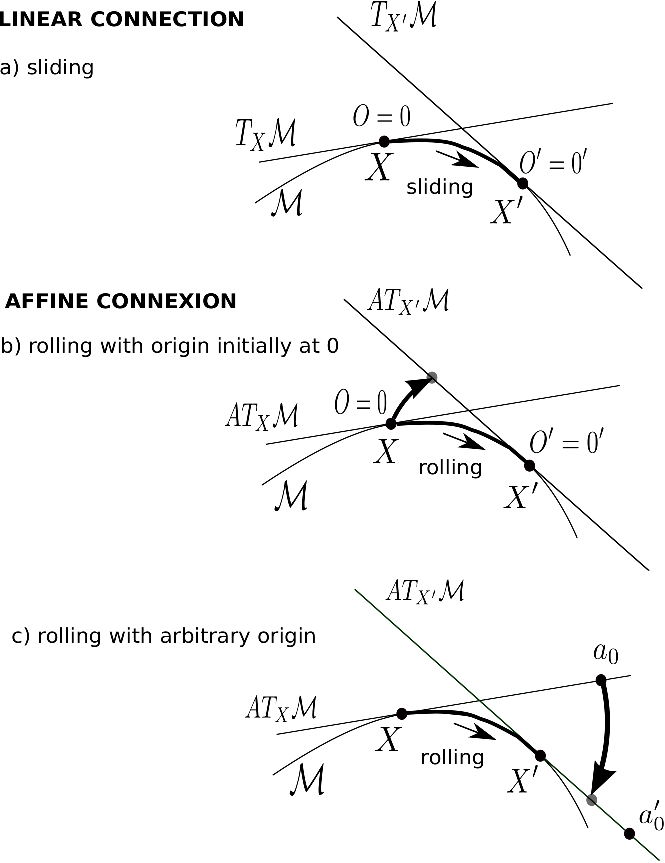} 
\caption{Linear and affine connection}
\label{fig ch1Figure1full}
\end{center}
\end{figure}

\section{Affine tensor analysis and principle of divergence free torsor}
\label{Section Affine tensor analysis and principle of divergence free torsor}

If there is a Lie group $G$, there is a $G$-structure and an Ehresmann connection \cite{Kobayashi}. For the affine group, it is an affine connection, according to the original definition in \'Elie Cartan's seminal paper \cite{Cartan 1923, Cartan 1924a}. First, we explain intuitively the idea on Figure \ref{fig ch1Figure1full}. A connection is a tool "to stick" the tangent spaces at two infinitesimally nearby points $\bm{X}$ and $\bm{X}'$. 
To construct a classical, linear connection $\nabla$, we slide $T_{\bm{X}}\mathcal{M}$ on the manifold to stick to $T_{\bm{X}'}\mathcal{M}$,  in such way that $\bm{0}$ coincides with $\bm{0}'$ (Figure \ref{fig ch1Figure1full}.a).

To construct an affine connection $\tilde{\nabla}$ with the origins at $\bm{0}$ and $\bm{0}'$, we roll the affine tangent space $AT_{\bm{X}}\mathcal{M}$ to stick to $AT_{\bm{X}'}\mathcal{M}$ but $\bm{0}$ does not coincide with $\bm{0}'$ anymore (Figure \ref{fig ch1Figure1full}.b). Working more generally with arbitrary origins $\bm{a}_0$ and $\bm{a}'_0$, the two origins do not coincide when we stick both affine tangent spaces (Figure \ref{fig ch1Figure1full}.c).

Then to define a covariant derivative $\tilde{\nabla}$ of an affine tensor, we consider an infinitesimal Galilean transformation where

\begin{itemize}
    \item the connection matrix 
    $$ \Gamma = (\Gamma^{\alpha}_{\mu\beta} dX^\mu)
    $$
    represents the infinitesimal motion of the basis and is physically identified to the gravitation, as in Relativity.
    \item the connection column vector
    $$ \Gamma_A = (\Gamma^{\alpha}_{A\mu} dX^\mu)
    $$
    represents the infinitesimal motion of the origin, physically the motion of the observer. Owing to (\ref{affine frame change}) and (\ref{C' = & P = & C =}), $C$ represents a vector field in a local chart, then it is meaningful to consider its covariant derivative $\nabla_{dX} C = dC + \Gamma(dX) \, C$. Then it can be proved that 
    \begin{eqnarray}
         \Gamma_A (dX) = dX -  \nabla_{dX} C
    \label{Gamma_A (dX) = dX - nabla_dX C}
    \end{eqnarray}
    in which the term $dX$ results from the rolling in the middle figure while the second term containing $C$ has to be added for the general case represented in Figure \ref{fig ch1Figure1full}.c. 
\end{itemize}
For a more rigorous presentation, the reader is referred for instance to (\cite{AffineMechBook}, Section 14.5.3).

On this ground, we define the covariant divergence of the torsor of a Cosserat medium. It is a field on the matter manifold $\mathcal{N}$, embedded into $\mathcal{M}$. Then the connection on $\mathcal{N}$ is in fact the pullback connection of $\tilde{\nabla}$. The affine covariant divergence of a vector-valued torsor $\bm{\tau}$ is a real-valued torsor, then decomposed as follow
\begin{eqnarray}
          \tilde{\bm{Div}}\,\bm{\tau} = \, _\gamma \tilde{\nabla}\,^\gamma T^\beta\,
         ( \bm{a}_0 \otimes \vec{\bm{e}}_\beta - \vec{\bm{e}}_\beta  \otimes\bm{a}_0) 
  + \, _\gamma \tilde{\nabla} \,^\gamma J^{\alpha\beta}\,\vec{\bm{e}}_\alpha  \otimes \vec{\bm{e}}_\beta
    \label{covariant divergence of a torsor}
    \end{eqnarray}
where the expression in local charts of the divergence of $T$ and $J$ components
\begin{eqnarray}
      _\gamma \tilde{\nabla}\,^\gamma T^\beta = 
               \frac{\partial(\, ^\gamma T^{\beta})}{\partial (\, ^\gamma \xi)} 
               +\, ^{\gamma}_{\gamma\rho} \Gamma\, ^\rho T^\beta 
               +\, ^\gamma T^\rho \, _\gamma U^\sigma \, \Gamma^\beta_{\sigma\rho}     
    \label{div T expression}
\end{eqnarray}
$$   _\gamma \tilde{\nabla}\;^\gamma J^{\alpha\beta}  =  
                                           \frac{\partial(\, ^\gamma J^{\alpha\beta} )}{\partial(\, ^\gamma \xi)}
                                         + \,^\gamma J^{\rho\beta}\,  _\gamma U^\sigma \, \Gamma^\alpha_{\sigma\rho}
                                         + \,^\gamma J^{\alpha\rho}\, _\gamma U^\sigma \, \Gamma^\beta_{\sigma\rho}
                                         + \,^\gamma_{\gamma\rho} \Gamma\, ^\rho J^{\alpha\beta}
$$
\begin{eqnarray}
     + \, _\gamma U^\sigma \, \Gamma^\alpha_{A\sigma}\,^\gamma T^\beta 
                                         - \,^\gamma T^\alpha \, _\gamma U^\sigma\, \Gamma^\beta_{A\sigma}   
    \label{div J expression}
\end{eqnarray} 
are obtained by standard techniques of tensor calculus. 


There are many Christoffel symbols but a lot of them vanish in classical mechanics. Indeed, for the Galilean gravitation, of the $40$ independent Christoffels, only $6$ remains
\begin{itemize}
    \item the 3 components of the gravity
    \begin{eqnarray}
         g^{i} = - \Gamma^{i}_{00}
    \label{g^i =}
    \end{eqnarray}
    \item and the 3 components of the axial vector $\Omega$ associated to the spinning, responsible of Coriolis' force
    \begin{eqnarray}
          \Omega^{i}_{j} = - \Omega^{j}_{i} = \Gamma^{i}_{0j} 
    \label{Omega^i_j =}
    \end{eqnarray}
\end{itemize}


We are able now to modelize the dynamics of Cosserat media compatible with Galilei principle of Relativity. In modern language, the balance equations must be covariant. Then, we claim the \textbf{principle of divergence free torsor:}  

\begin{eqnarray}
    \tilde{\bm{Div}}\,\bm{\tau}  = \bm{0}
\label{principle of divergence free torsor}    
\end{eqnarray}

\vspace{0.1cm}

It is abstract but covers a broad spectrum of applications. As the torsor has 10 independent components, we deduce from this general principle a set of 10 scalar balance equations, whatever the dimension of the matter manifold. In fact, for each kind of medium, we obtain family of balance equation sets, depending on the choice of the local charts and origins of the affine frames. In particular, we shall indicate whether we are working in a proper frame or not. 

\section{Dynamics of a bulky body: the case of a Cauchy medium}
\label{Section Dynamics of a bulky body: the case of a Cauchy medium}

For easiness, we use once again coordinates (\ref{xi^alpha = X^alpha}) such that the embedding $i:\mathcal{N} \rightarrow \mathcal{M}$ is represented by the identity and we use the index convention (\ref{convention index at the right}). Then the tangent map $U$ to $i$ has components $U^\alpha_\beta = \delta^\alpha_\beta$. Hence the expressions of the components of the torsor divergence (\ref{div T expression}) and (\ref{div J expression}) are simplified
\begin{equation}
   \tilde{\nabla}_\gamma T^{\beta\gamma} = 
               \frac{\partial T^{\beta\gamma}}{\partial X^\gamma} 
               + \Gamma^\beta_{\gamma\rho} T^{\rho\gamma} 
               + \Gamma^\gamma_{\gamma\rho} T^{\beta\rho} \ ,
\label{tilde(div) T^(alpha gamma)} 
\end{equation}
\begin{equation}
   \tilde{\nabla}_\gamma J^{\alpha\beta\gamma} = \frac{\partial J^{\alpha\beta\gamma} }{\partial X^\gamma}
                                         + \Gamma^\alpha_{\gamma\rho} J^{\rho\beta\gamma} 
                                         + \Gamma^\beta_{\gamma\rho} J^{\alpha\rho\gamma}
                                         + \Gamma^\gamma_{\gamma\rho} J^{\alpha\beta\rho}
                                         + \Gamma^\alpha_{A\gamma} T^{\beta\gamma} 
                                         - T^{\alpha\gamma} \Gamma^\beta_{A\gamma}\ .
\label{tilde(div) J^(alpha beta gamma)} 
\end{equation}

Applying the principle of divergence free torsor, we obtain the 10 balance equations in a covariant form, \textit{i.e.} compatible with Galilei principle of relativity
\begin{itemize}
    \item for the $T$ components, Euler's equations of fluids:
    \begin{itemize}
        \item  the balance of mass
        $$ \dfrac{\partial \rho}{\partial t} + div\,(\rho\,v) = 0
        $$
        \item the balance of linear momentum
        $$ \rho\ \left[ \dfrac{\partial v}{\partial t} +  \dfrac{\partial v}{\partial x}\,v 
          \right] =  div\ \sigma + \rho (g - 2\ \Omega \times v)
        $$
        with at the left the Eulerian acceleration containing the advection term, and right Coriolis' acceleration that appear naturally. This simplified expression is obtained by using the balance of mass (\cite{AffineMechBook}, Section 5.3.3).
    \end{itemize}
    \item for the $J$ components, the symmetry of the stress-mass tensor of the Cauchy media:
    $$ T^{\alpha\beta} = T^{\beta\alpha} 
    $$
    by working in a proper frame, then $J^{\alpha\beta\gamma} = 0$,  $C^\alpha = 0$ and, owing to (\ref{Gamma_A (dX) = dX - nabla_dX C}),  $\Gamma^\alpha_{A\beta} = \delta^\alpha_\beta$ .
\end{itemize}

\section{Dynamics of a pointwise object}
\label{Section Dynamics of a pointwise object}

In this special case, the torsor is scalar valued, the left indices and some terms disappeared of the expressions. Besides, we can use for easiness matrix notations as in Section \ref{Subsection Transformation laws of affine tensors}. Then the components of the covariant derivative of the torsor read
$$  \tilde{\nabla}_{dX} T = dT + \Gamma (dX)\,T,\qquad 
   \tilde{\nabla}_{dX} J = dJ + \Gamma(dX)\,J + J\,(\Gamma(dX))^T +\Gamma_A(dX)\,T^T - T\,(\Gamma_A(dX))^T\ 
$$
For this kind of medium, as in Section \ref{Subsection The torsor of a pointwise object}, we do not work in a proper frame. Owing to (\ref{tilde(P) = ((1 0)(C P)) =}), we consider the translation
\begin{eqnarray}
     C =\left( \begin{array} {c}
                       0 \\
                       x \\
             \end {array} \right) 
\label{C = (0 x)}
\end{eqnarray}
Taking into account the Galilean gravitation (\ref{g^i =}, \ref{Omega^i_j =}), (\ref{Gamma_A (dX) = dX - nabla_dX C}) gives
\begin{eqnarray}
     \Gamma_A = d \left( \begin{array} {c}
                       t               \\
                       x \\
                    \end {array} \right)
             - d\left( \begin{array} {c}
                       0               \\
                       x \\
                    \end {array} \right)
              - \left( {{\begin{array}{cc}
                         0 &  0               \\
                     \Omega \times dx - g\,dt  &  j(\Omega)\,dt   \\
                      \end{array} }} \right)\,
              \left( \begin{array} {c}
                       0               \\
                       x \\
                    \end {array} \right)
             = \left( \begin{array} {c}
                         dt                               \\
                       - \Omega\times x\,dt \\
                    \end {array} \right)
\label{Galilean Gamma_A}
\end{eqnarray}
where $j(\Omega)$ is the skew-symmetric matrix of components (\ref{Omega^i_j =}). From the principle of divergence free torsor (\ref{principle of divergence free torsor}), we get also 10 balance equations in a covariant form, \textit{i.e.} with Coriolis' acceleration:
\begin{itemize}
    \item the balance of mass
    $$ \dot{m} = 0
    $$
    \item the balance of linear momentum
    $$ \dot{p}  =  m \, (g - 2\ \Omega \times v) 
    $$
    This covariant equation with Coriolis force was proposed by J.-M. Souriau (\cite{SSD, SSDEng}, Equations (12.44) and (12.47)).
    \item the balance of position quantity
    \begin{eqnarray}
        \dot{q} = p
    \label{dot(q) = p}
    \end{eqnarray}
    This equation is obvious since $m$ is time-independent then the time-derivative of $m x$ is $m v$
    \item the balance of angular momentum
    $$ \dot{l} + \Omega \times l_0 
    = x \times  m \, (g - 2\ \Omega \times v) 
    $$
    This equation is not trivial. It allows to understand the motion of  the Lagrange top (\cite{AffineMechBook}, Sections 5.5). This simplified expression is obtained by using (\ref{Galilean Gamma_A}) and Jacobi's identity for the cross-product (\cite{de Saxce 2011} and \cite{AffineMechBook}, Section 5.3.3). If we rather had worked in a proper frame, we would have obtained the motion \textit{\`a la Poinsot}, for instance to describe the attitude of a satellite (\cite{AffineMechBook}, Sections 5.4)
\end{itemize}


\section{Dynamics of a 1D Cosserat medium}
\label{Section Dynamics of a 1D Cosserat medium}

Once again, we apply the principle of divergence free torsor (\ref{principle of divergence free torsor}), using (\ref{div T expression}) and (\ref{div J expression}), that leads to equations that are valid not only for solids (beams, arm, string) but also for fluids (flows in a pile, water jet):
\begin{itemize}
    \item The balance of mass 
    $$ \dfrac{\partial \rho_l}{\partial t} + \dfrac{\partial }{\partial s} (\rho_l v_t) = 0
    $$
    is a simple confirmation of a condition that can be obtained by straightforward methods
    \item In the balance of linear momentum, 
    $$ \rho_l \left[ \dfrac{\partial v}{\partial t} +  \dfrac{\partial v}{\partial s}\,v_t 
          \right]
                   = \dfrac{\partial F}{\partial s} 
                     + \rho_l (g - 2\ \Omega \times v) 
    $$
    we recognize at the left the advection term in the direction tangent to the body and right the Coriolis' effect, then in a covariant form
    \item The balance of position quantity
    \begin{eqnarray}
         \dfrac{\partial q}{\partial t} + \dfrac{\partial l_\star}{\partial s} = \rho_l v
    \label{1D - balance of position quantity}
    \end{eqnarray}
    It was obtain with the motion of the origin given by (\ref{C = (0 x)}) and the corresponding connection column vector (\ref{Galilean Gamma_A}) that reads in tensor notations: $ \Gamma^0_{A0} = 1, \Gamma^i_{A0} =- \Omega^i_j x^j, \Gamma^\alpha_{Aj} = 0$, as in the case of the pointwise object for which the equation (\ref{dot(q) = p}) was obvious. For the 1D body, it is non trivial and useful. Indeed, the 3D-to-1D reduction is only an approximation. At the left, $q$ contains the mass-center, that is the average position on the cross-section. The problem is that the time-derivative of the average of the position is not the average of the time-derivative, that is the velocity at right (for instance a flow in an elbow). This balance equation is valuable because it allows to restore the compatibility destroyed by the averaging.
    \item The balance of angular momentum
    \begin{eqnarray}
         \dfrac{\partial l}{\partial t}
    + \Omega \times l 
    + l_\star \times (\Omega \times n) 
    = - \dfrac{\partial M_\star}{\partial s}
    + n \times F 
    \label{1D - balance of angular momentum}
    \end{eqnarray}
    obtained in a proper frame as usual in beam theory. In the limit case of the Statics, the left hand side disappears and we recover the classical relation between bending moments and shear forces of classical beam theory. At the left, the two former terms are similar to those of the dynamics of the top. The last one containing the tangent angular momentum is absolutely non trivial. 
\end{itemize}

\section{Dynamics of 2D Cosserat media}
\label{Section Dynamics of 2D Cosserat media}

The dynamics of 2D Cosserat media is more complex and is developed in \cite{de Saxce 2003}. We give only the final results of tricky calculations. The reader is referred to \cite{de Saxce 2003} for technical details. The modelling of the matter motion is similar to the one of 1D Cosserat media with the additional difficulty of working as usual with curvilinear coordinates adapted to the shape of the surface. In the limit case of statics, we recover the standard system of 10 equilibrium equations \cite{Green Zerna 1968}. As for the 1D Cosserat medium, the interest of the present geometric formulation is elsewhere, in the balance equations of the dynamics where appear non trivial new terms that could be useful for the modelling of impact problems in solid mechanics and for fluids. 

The local representation of the embedding of the thin 3D body (called shell for simplicity even if it is a fluid) into the space-time is supposed to be of the following form:
 $$ i: \bar{\xi} =  \left( \begin{array} {c}
        t    \\
        \bar{s}    \\
\end {array} \right) \mapsto \bar{X} = 
 \left( \begin{array} {c}
        t    \\
        \bar{x}    \\
\end {array} \right)  =
 \left( \begin{array} {c}
        t    \\
        \psi (t, \bar{s})   \\
\end {array} \right) =
 \left( \begin{array} {c}
        t    \\
        x (t) + R(t) \, \bar{s}   \\
\end {array} \right) 
 $$
where the point on the middle surface $x$, identified to the 2D body which idealizes the 3D one, and the rotation $R$ describe the rigid motion of the shell element around this point. Its Poisson's vector $\varpi$ is such that $\dot{R} = j (\varpi) \, R$. Concerning the Latin and Greek indices, we use the same convention  as in Section \ref{Subsection 3D-to-1D reduction}, with the additional rule that $a, b, c \in \left\lbrace 1, 2 \right\rbrace $. A classical tool of the theory of surface is the tangent plane at the current point $x$. In order to separate in-plane and off-plane components of the torsor and the momentum balance, we introduce another coordinate systems in which the embedding is represented by
$$ i: \bar{\xi}' =  \left( \begin{array} {c}
        t    \\
        \theta^1    \\
        \theta^2    \\
\end {array} \right) \mapsto 
\bar{X}' = 
 \left( \begin{array} {c}
        t    \\
        \theta^1    \\
        \theta^2    \\
        \theta^3    \\
\end {array} \right) 
$$
where $\theta^i$ are curvilinear coordinates adapted to the body in such way that $(\theta^a)$ parameterize the middle surface of equation $\theta^3 = 0$ and the tangent map $U$ to the embedding $i$ is very simple
$$ _0 U^0 = 1, \quad
   _a U^b = \delta^b_a, \quad
    _a U^0 = \, _a U^3 = \, _0 U^i = 0
$$
They are related to the former coordinates by
$$ \bar{x}^i = x^i (t, \theta^a) + \theta^3 n^i (t, \theta^a)
$$
For convenience, following \cite{Naghdi 1972}, we choose $n^i$ such that 
$$  \pi^i_a \delta_{ij} n^j = 0, \qquad
    n^i \delta_{ij} n^j = 1
$$
with
$$ \pi^i_a = \frac{\partial x^i}{\partial \theta_a}, \qquad
v^i = \frac{\partial x^i}{\partial t}, \qquad
w^i = \frac{\partial n^i}{\partial t} = (\varpi \times n)^i
$$
then $n$ represents the unit vector normal to the tangent plane and $\pi$ the projector onto the tangent plane. In the adapted coordinates and assuming that the 3D body is a Cauchy medium, the $T$ components of the torsor are
$$ \bar{T}^{ab} = \rho \, (\theta^3)^2 w^a \, w^b 
                                    - \bar{\sigma}^{ab}, \qquad
   \bar{T}^{a0} = \bar{T}^{0a} = \rho \, \theta^3 w^a 
$$
$$ \bar{T}^{a3} = \bar{T}^{3a} =  - \sigma^{a3}, \qquad
   \bar{T}^{30} = \bar{T}^{03} = 0, \qquad
   \bar{T}^{00} = \rho
$$

By analogy with the Section \ref{Subsection 3D-to-1D reduction}, we perform a \textbf{3D-to-2D reduction} that, by projection and integration over the thickness $h$ leads to the components of the 2D body torsor
$$ ^a T ^b = \frac{\rho \, h^3}{12}\,  w^a \, w^b 
                                    - N^{ab}, \qquad
   ^a T^3 = - Q^a, \qquad 
   ^0 T^0 = \rho_s, \qquad
    ^a J^{b3} = - ^a J^{3b} = M^{ab}
$$
where we define the shell variables
\begin{itemize}
    \item the mass density $\rho_s = \rho \, h$
    \item the membrane force density $N^{ab} 
     = \int^{h/2}_{- h/2} \bar{\sigma}^{ab} \, \mbox{d} \theta^3$ 
    \item the shear force density $Q^a 
     = \int^{h/2}_{- h/2} \bar{\sigma}^{a3} \, \mbox{d} \theta^3$ 
    \item the bending/torsion moment density $M^{ab} 
     = \int^{h/2}_{- h/2} \theta^3\bar{\sigma}^{ab} \, \mbox{d} \theta^3$ 
\end{itemize}
where all the densities are per surface unit. 

Introducing the first fundamental form  of the surface $a = \pi \, \pi^T$, $c = a^{-1} \pi$, the in-plane component of the gravity $g_t = c \, g$ and the off-plane component $g^3 = n^T g$, the symbol $\varepsilon_{cb}$ such that $\varepsilon_{12} = - \varepsilon_{21} = 1, \varepsilon_{11} =  \varepsilon_{22} = 0$ , the Christoffel's symbols read:
$$ \Gamma^a_{00} = - (g_t)^a, \quad 
  \Gamma^3_{00} = - g^3, \quad
\Gamma^a_{bc} = c^a_i \frac{\partial \pi^i_b}{\partial \theta^c}, \quad
\Gamma^3_{ab} = n^j \delta_{ji} \frac{\partial \pi^i_b}{\partial \theta^a} = b_{ab}
$$
where we recognize the components $b_{ab}$ of the second fundamental form, and:
$$ \Gamma^a_{b0} = c^a_i \left( \frac{\partial \pi^i_b}{\partial t} + \Omega^i_j \pi^j_b\right) = \Phi^a_b, \quad 
\Gamma^a_{30} = c^a_i \left( \frac{\partial n^i}{\partial t} + \Omega^i_j n^j\right) = \Phi^a
$$
$$ \Gamma^3_{b0} = n^j \delta_{ji} \left( \frac{\partial \pi^i_b}{\partial t} + \Omega^i_j \pi^j_b\right) = \Phi_b
$$
In the sequel, we use standard shortcuts 
$$ N^{ba}  \mid_b = \frac{\partial N^{ba} }{\partial \theta^b}
 + \Gamma^a_{bc} N^{bc} + \Gamma^c_{cb} N^{ba} , \quad
 Q^b  \mid_b = \frac{\partial Q^b }{\partial \theta^b}
 + \Gamma^c_{bc} Q^b
$$
The balance equations of angular momentum are obtained in a proper frame as usual in classical shell theory. Applying the principle of divergence free torsor (\ref{principle of divergence free torsor}), we obtain:
\begin{itemize}
    \item the balance of mass
    $$ _\gamma \tilde{\nabla} ^\gamma T^0 = \frac{\partial \rho_s}{\partial t} + \Phi^a_a \rho_s = 0
    $$
    \item the balance of in-plane linear momentum
    $$ _\gamma \tilde{\nabla} ^\gamma T^a = 
    \left( N^{ba} - \frac{\rho \, h^3}{12}\,  w^a \, w^b  \right) \mid_b - b^a_b Q^b +  c^a_i \left( \rho_s \,(g^i - 2 \, \Omega^i_k v^k) - \rho_s \frac{\partial v^i}{\partial t} \right)  = 0
    $$
     \item the balance of off-plane linear momentum
    $$ _\gamma \tilde{\nabla} ^\gamma T^3 = b_{ab} 
    \left( N^{ba} - \frac{\rho \, h^3}{12}\,  w^a \, w^b  \right) + Q^b \mid_b + n^j \delta_{ij} \left(\rho_s  (g^i - 2 \, \Omega^i_k v^k) - \rho_s \frac{\partial v^i}{\partial t} \right) = 0
    $$
    \item the balance of in-plane angular momentum
    \begin{eqnarray}
        _\gamma \tilde{\nabla} ^\gamma J^{21} = \varepsilon_{cb} 
    \left( N^{cb} - b^c_a M^{ab}  - \frac{\rho \, h^3}{12}\,  (\Phi^c + w^c)\, w^b \right) = 0 
    \label{2D - balance of in-plane angular momentum}
    \end{eqnarray}
    that generalize to the dynamics the symmetry relation of the usual shell theory \cite{Valid 1995}
    \item the balance of off-plane angular momentum
    \begin{eqnarray}
         _\gamma \tilde{\nabla} ^\gamma J^{b3} = M^{ab} \mid_b - Q^b
      - \frac{\rho \, h^3}{12}\, \left( \frac{\partial w^b}{\partial t}+ \Phi^b_a w^a + \Phi^c_c w^b) \right) = 0
    \label{2D - balance of off-plane angular momentum}
    \end{eqnarray}
\end{itemize}
Finally, it can be seen that the remaining equations are automatically satisfied:
$$ _\gamma \tilde{\nabla} ^\gamma J^{03} = 0, \qquad
   _\gamma \tilde{\nabla} ^\gamma J^{b0} = 0
$$

\section{Dynamics of a 3D Cosserat medium}
\label{Section Dynamics of a 3D Cosserat medium}

While the relevance of Cosserat theory for slender and thin bodies is well-established, its application to 3D bodies causes problems concerning the kind of boundary conditions to prescribe and the measure of the material parameters. Nonetheless, starting from the decade of 50s, there is a significant literature devoted to to its use for the Statics of solid bodies.  There are often qualified of micropolar materials. The purposes were to model multilayer materials \cite{Forest 1998}, strongly heterogeneous materials \cite{Neff 2007}, material with random structures \cite{Forest 2000}, granular materials \cite{Kruyt 2003, Ching Chang 2005}, diffuse interfaces \cite{Ask 2018} and shear bands. Another attempt is applying Cosserat theory to 3D fluids to model the turbulence thanks to the introduction, besides of the stress tensor, of a stress moment tensor and corresponding constitutive equations deriving from a dissipation potential for the viscous effects \cite{Eringen 1964, Eringen 2003}.

Our idea is that the range of application is much more open if we extend this Cosserat framework to the Dynamics. It is what we would like to explore thanks to the principle of divergence free torsor. To reduce the number of indices, we use engineer notations analogous to the ones of the 1D Cosserat media
$$ q^i =\, J^{i00},\quad 
\quad l^i = J^{kl0}, \quad
l^{ir}_\star = J^{i0r}, \quad
\quad M^{ir}_\star = J^{klr} 
$$
where $(ikl)$ is a cyclic permutation of $(123)$. The stress moment tensor is given by the components $M^{ir}_\star$ but their are other components $q^i, l^i$ and $l^{ir}_\star$ which enhance the model of the static Cosserat model. 
We apply the principle of divergence free torsor (\ref{principle of divergence free torsor}) and calculate the divergence of the torsor in a proper frame, using (\ref{tilde(div) T^(alpha gamma)}) and (\ref{tilde(div) J^(alpha beta gamma)}), that leads to the following equations with:
\begin{itemize}
        \item  the balance of mass
        $$ \dfrac{\partial T^{00}}{\partial t} + 
        \dfrac{\partial T^{0i} }{\partial x^i} = 0
        $$
        \item the balance of linear momentum
        $$ \dfrac{\partial T^{i0}}{\partial t} + 
        \dfrac{\partial T^{ij} }{\partial x^j} 
        - (T^{00} g^i - \Omega^i_j (T^{0j} + T^{j0})) = 0
        $$
        \item the balance of position quantity
        $$ 
        \dfrac{\partial q^i}{\partial t} + \Omega^i_j q^j
        + \dfrac{\partial l^{ir}_\star}{\partial x^r} 
        + T^{0i} - T^{i0} = 0
        $$
        \item the balance of angular momentum
        $$
        \dfrac{\partial l^k}{\partial t} 
        + \dfrac{\partial M^{km}_\star}{\partial x^m} 
        - (q \times g)^k
        + \Omega^j_p l^q - \Omega^i_r l^s
        + \Omega^j_r l^{ir}_\star - \Omega^i_r l^{jr}_\star 
        + T^{ji} - T^{ij} =0 $$
        where $(ijk)$, $(ipq)$ and $(jrs)$ are cyclic permutations of $(123)$ and summations have to be performed on all triples $(ipq)$ and $(jrs)$.
    \end{itemize}


\section{Conclusions and perspectives}
\label{Section Conclusions and perspectives}

In this paper, we introduced a geometric framework based on the concept of affine tensors. In the context of Dynamics, we gave a survey of applications covering a large range from the pointwise object to the 3D body. The essential point is the consideration of the time as as a coordinate in its own right, on the same footing as the three spatial coordinates. 
Particular quantities (mass, force, stress, moments, dynamic momenta and so on...) that are introduced in the course of standard textbooks "drop by drop", are gathered here in big tensors of higher dimensions and ranks. These quantities are working together as a team.

On this ground, we proposed a general principle of divergence free torsor that ensure the consistency with Galilei principle of relativity. For each kind of application, it provides a set of ten balance equations of which many are well-known but some of them are new, namely for the pointwise object (Section \ref{Section Dynamics of a pointwise object}), Equations  (\ref{1D - balance of position quantity}) and (\ref{1D - balance of angular momentum}) for the 1D Cosserat medium (Section \ref{Section Dynamics of a 1D Cosserat medium}), Equations (\ref{2D - balance of in-plane angular momentum}) and (\ref{2D - balance of off-plane angular momentum}) for the 2D Cosserat medium (Section \ref{Section Dynamics of 2D Cosserat media}), and the balance equations for the 3D Cosserat medium (Section \ref{Section Dynamics of a 3D Cosserat medium}). Note also the simple Equations (\ref{T = int_S Pi bar(T) dS}) and (\ref{J = int_S Pi bar(J) dS}) of the 3D-to-1D reduction which are original. Take car also that the statical force (\ref{F = int_S (sigma n - rho (v_t - bar(v)_t) (v - bar(v)) dS}) contains a dynamical contribution. Finally, even if we considered only applications to the classical mechanics, the general framework proposed in this paper is also valid in General relativity, in particular Equations (\ref{covariant divergence of a torsor}), (\ref{div T expression}), (\ref{div J expression}) and the principle (\ref{principle of divergence free torsor}).

\vspace{0.3cm}

As perspectives, we would like to discuss the classification of affine tensors that we consider important for the Mechanics.
The topics of torsors was deeply developed in this paper 
but it remains other kinds of affine tensors which deserve to be studied in the future.
\begin{itemize}
    \item \textbf{Torsors} (2-contravariant) : 
    $$ (\bm{\Psi}, \hat{\bm{\Psi}}) \mapsto \bm{\tau} (\bm{\Psi}, \hat{\bm{\Psi}} )\in T_\xi \mathcal{N}
    $$
    They are skew-symmetric.
    \vspace{0.2cm}
    \item \textbf{Co-torsors} (2-covariant) : 
$$ (\bm{a}, \hat{\bm{a}} ) \mapsto \bm{\gamma} (\bm{a}, \hat{\bm{a}} )\in T^*_\xi \mathcal{N}
$$
Their two arguments are points and the co-torsors are bi-affine and skew-symmetric. They allow to describe the kinematics, for instance the one of a rigid body, and are in duality with torsors (see \cite{AffineMechBook}, Section 5.1.3., page 76-80). 
We can combine the method of virtual powers with the affine tensors by defining the gradient of a co-torsor.
    \vspace{0.2cm}
    \item \textbf{Momentum tensors} (1-covariant and 1-contravariant)  :
$$ (\vec{\bm{V}}, \bm{\Psi}) \mapsto \bm{\mu} (\vec{\bm{V}}, \bm{\Psi}))\in T_\xi \mathcal{N}
$$
Their 2 arguments are a vector and an affine form. They are bilinear but not skew-symmetric. They can be decomposed as follows:
$$ \bm{\mu} = \,^\gamma \bm{\mu}\,_\gamma\vec{\bm{\eta}},\qquad 
^\gamma \bm{\mu} = \bm{e}^\beta \otimes 
        (\,^\gamma \Sigma_\beta \, \bm{a}_0 + \,^\gamma M^\alpha_\beta \,  \vec{\bm{e}}_\alpha)
$$
Their components can be stored in a row $^\gamma \Sigma$ and a square matrix $^\gamma M$. 
We can recognize in the couple $(^\gamma \Sigma, ^\gamma M)$ an element of the dual of the Lie algebra of the affine group. It is remarkable to observe that their transformation law is just the coadjoint representation of the affine group $\mathbb{GA} (m) $ where $m = \mbox{dim}(\mathcal{M})$. Then the system of components of $\bm{\mu}$ can be considered as the value of the momentum map in symplectic mechanics, then the name of momentum tensors.
    \vspace{0.2cm}
    \item \textbf{Strain tensors} (1-contravariant  and 1-covariant)  : 
    $$ (\bm{\Phi}, \bm{a} ) \mapsto \bm{\omega} (\bm{\Phi}, \bm{a} ) \in T^*_\xi \mathcal{N} 
    $$
    Their arguments are a linear form and a point. They are linear and affine. They represent generalized deformations of Cosserat media in duality with momentum tensors.    
\end{itemize}
In the future we hope to develop considerations 
relative to these kind of tensors. The combination with a variational formulation as in  \cite{Boyer 2009, Boyer 2017} seems to be promising. 

\vspace{0.3cm}

\textbf{Acknowledgements}

\vspace{0.3cm}

This work were presented at the Workshop on "Geometrical Aspects of Material Modelling" (GAMM), 21-23 August 2024, Madrid (Spain). The author thanks the organizers for this worthwhile initiative and their kind invitation.


\end{document}